\documentclass[a4paper, 12pt]{amsart}
\usepackage{amsmath}

\usepackage{amssymb,latexsym}
\usepackage[dvips]{graphicx}
\usepackage[latin1]{inputenc}

\newcommand{\R}{\mathbb R}
\newcommand{\Sp}{\mathbb S}
\newcommand{\Z}{\mathbb Z}
\newcommand{\N}{\mathbb N}

\newtheorem{theorem}{Theorem} [section]
\newtheorem{lemma}{Lemma} [section]
\newtheorem{proposition}{Proposition} [section]
\newtheorem{corollary}{Corollary} [section]
\newtheorem{definition}{Definition} [section]

\newtheorem{remark}{Remark}[section]

\begin{document}

\title[Critical metrics of  the Laplacian eigenvalues]
{Laplacian eigenvalues functionals and metric deformations on compact manifolds}
\author{Ahmad El Soufi and Sa\"{\i}d Ilias}

\date{}

\address{Universit\'e de Tours, Laboratoire de Math\'ematiques
et Physique Th\'eorique, UMR-CNRS 6083, Parc de Grandmont, 37200
Tours, France} \email{elsoufi@univ-tours.fr, ilias@univ-tours.fr}


\begin{abstract}
In this paper, we investigate critical points of the Laplacian's eigenvalues considered as functionals on the space of Riemmannian metrics or a conformal class of metrics on a compact manifold.
We obtain necessary and sufficient conditions for a metric to be a critical point of such a functional. We derive specific   consequences concerning possible locally maximizing metrics. We also characterize critical metrics of the ratio of two consecutive eigenvalues.

\end{abstract}

\maketitle


\section {Introduction}
The eigenvalues of the Laplace-Beltrami operator associated with a Riemannian metric on a closed manifold are among the most natural global Riemannian invariants defined independently from curvature. One of the main topics in Spectral Geometry is the study of uniform boundedness of eigenvalues under some constraints and the finding out of eventual extremal metrics. Let us start by recalling some important results in this direction, where the eigenvalues are considered as functionals on the set of Riemannian metrics of fixed volume.

In all the sequel, we will denote by $M$ a compact smooth manifold of dimension $n \geq 2$ and, for any Riemannian metric $g$ on $M$, by $$ 0<\lambda_{1}(g) \leq
\cdots \leq \lambda_{k}(g) \leq \cdots\rightarrow \infty$$ the sequence of eigenvalues of the Laplacian $\Delta_{g}$ associated with $g$, repeated according to their multiplicities. Notice that $\lambda_{k}$ is not invariant under scaling (i.e., $\forall c>0, \, \lambda_{k}(cg)= {1 \over c} \lambda_{k}(g)$). Hence, a normalization is needed and it is common to restrict the functional $\lambda_{k}$ to the set $\mathcal R (M)$ of Riemannian metrics of fixed volume. We will denote by $C(g)$ the set of Riemannian metrics conformal to $g$ and having the same volume as $g$.

The result that, in dimension $2$, $\forall k \geq 1$, the functional $\lambda_{k}$ is uniformly bounded on the set of metrics of fixed area is due to Korevaar \cite{K}, after been proved for $k=1$, by Hersch \cite{H} in the case of the 2-sphere $\Sp^2$, and by Yang and Yau \cite {YY} and Li and Yau \cite {LY} for all compact surfaces (see also \cite{EI} for an improvement of the Yang-Yau upper bound of $\lambda_{1}$ in terms of the genus).

The situation differs in dimension $n \geq 3$. Indeed, based on earlier results on spheres obtained by many authors \cite{TA,M,BB,U,MOU}, Colbois and Dodziuk \cite{CD} proved that, for any compact manifold $M$ of dimension $n \geq 3$, the functional $\lambda_{1}$ is unbounded on $\mathcal R (M)$.

However, $\forall k \geq 1$, the functional $\lambda_{k}$ becomes uniformly bounded when restricted to a conformal class of metrics of fixed volume $C(g)$. This result was first proved for $k=1$  by the authors \cite{EI0} (see also \cite{FN}) and, for any $k \geq 1$, by Korevaar \cite{K} (see also \cite{GY}).

Existence results for maximizing metrics are available in only few situations and concern exclusively the first eigenvalue functional. Hersch \cite{H} proved that the standard metric is the only maximizing metric for $\lambda_{1}$ on the $2$-sphere $\mathbb{S}^{2}$. The same result holds for the standard metric of the real projective plane $\R P^{2}$ (Li and Yau \cite {LY}). Nadirashvili \cite{N1} outlined a proof for the existence of maximizing metrics for $\lambda_{1}$ on the $2$-Torus and the Klein bottle (see \cite{G1,G2} for additional details on Nadirashvili's paper).

In higher dimensions, the authors  \cite{EI0} gave a sufficient condition for a metric $g$ on $M$ to maximize  $\lambda_{1}$ in its the conformal class  $C(g)$. This condition is fulfilled in particular by the standard metric of the sphere $\mathbb{S}^{n}$ (which enabled us to answer Berger's problem concerning the maximization of $\lambda_{1}$ restricted to the standard conformal class of $\mathbb{S}^{n}$), and more generally, by the standard metric of any compact rank-one symmetric space. The flat metrics $g_{sq}$ and $g_{eq}$ on the 2-Torus $\mathbb{T}^{2}$ associated with the square lattice $\Z^{2}$ and the equilateral lattice $\Z (1,0)\oplus \Z (1/2,\sqrt{3}/2)$ respectively are also maximizing metrics of $\lambda_{1}$ in their conformal classes (see \cite{EI0,LY}).\\

In this paper we address the following natural questions:
\begin{itemize}
	\item [(1)] \emph{What about critical points of the functional \;$g \mapsto \lambda_{k}(g)$} ?
	\item [(2)] \emph{How to deform a Riemannian metric $g$ in order to increase, or decrease, the $k$-th eigenvalue $\lambda_{k}$} ?
\end{itemize}

Despite the non-differentiability of the functional $\lambda_{k}$ with respect to metric deformations, perturbation theory enables us to prove that, for any analytic deformation $g_{t}$ of a metric $g$, the function $t\mapsto\lambda_{k}(g_{t})$ always admits left and right derivatives at $t=0$ (see Theorem \ref{main} (i) below). 
Moreover, these derivatives can be expressed in terms of the eigenvalues of the following quadratic form
$$ Q_{h}(u)=-\int_M\left( du\otimes du+ {1 \over 4} \Delta_{g}u^{2}g,h\right) v_g,$$
with $h={d \over dt}g_{t}\big|_{t=0}$, restricted to the eigenspace $E_k(g)$ (Theorem \ref{main} and Lemma \ref{lem 21}). 
This enables us to give partial answers to Question (2) above (Corollary \ref{cor 21}). In particular, if $Q_{h}$ is positive definite on $E_k(g)$, then $\lambda_{k}(g_{-t})<\lambda_{k}(g)<\lambda_{k}(g_{t})$ for all $t\in(0,\varepsilon)$, for some positive $\varepsilon$.

Concerning Question (1), the existence of left and right derivatives for $t\mapsto\lambda_{k}(g_{t})$, suggests the following natural notion of criticality. Indeed, a metric $g$ will be termed \emph{critical} for the functional $\lambda_{k}$ if, for any volume-preserving deformation $g_{t}$ of $g$, one has 
$$ {d \over dt}\lambda_{k}(g_{t})\big|_{t=0^-} \times {d \over dt}\lambda_{k}(g_{t})\big|_{t=0^+} \leq 0 \, ; $$
this means that either
$$\lambda_{k}(g_{t}) \leq \lambda_{k}(g) +o(t) \quad \hbox{or}\quad  \lambda_{k}(g_{t}) \geq \lambda_{k}(g) +o(t).$$
It is clear that if $g$ is a locally maximizing or a locally minimizing metric of $\lambda_{k}$, then $g$ is a critical metric for $\lambda_{k}$ in the previous sense. 

In an earlier work (\cite{EI1} and \cite{EI2}), we treated the particular case $k=1$ and gave a necessary condition ( \cite[Theorem 1.1]{EI1} and  \cite[Theorem 2.1]{EI2})
as well as a sufficient condition ( \cite[Proposition 1.1]{EI1} and  \cite[Theorem 2.2]{EI2}) for a metric $g$ to be critical for the functional $\lambda_{1}$ or for the restriction of the $\lambda_{1}$ to a conformal class.

In section 3 and 4 below,  we extend the results of \cite{EI1} and \cite{EI2} to higher order eigenvalues, and weaken the sufficient conditions given there. Actually, as we will see, in many cases, the necessary condition of criticality is also sufficient (Theorem \ref{th 31} and Theorem \ref{th 41}).
Given a metric $g$ on $M$, we will prove that
\begin{itemize}
\item \emph{Necessary conditions} (Theorem\, \ref {th 31}(i) and  Theorem \,\ref{th 41}(i)):
\begin{enumerate}
\item[i)] If $g$ is critical for the functional $\lambda_{k}$, then there exists a finite family $\left\{u_{1}, \cdots , u_{d}\right\}$
of eigenfunctions associated with $\lambda_{k}$ such that
$$ \sum_{i \leq d} du_{i} \otimes du_{i} =g.$$
\item[ii)] If $g$ is a critical metric of the functional $\lambda_{k}$ restricted to the conformal class of $g$, then there exists a finite family $\left\{u_{1}, \cdots , u_{d}\right\}$ of eigenfunctions associated with $\lambda_{k}$ such that $$\sum_{i \leq d} u_{i}^{2}=1. $$
\end{enumerate}
\item \emph{Sufficient conditions} (Theorem \ref {th 31}(ii) and Theorem \ref{th 41}(ii)): \\
 If $\lambda_{k}(g)>\lambda_{k-1}(g)$ or $\lambda_{k}(g)<\lambda_{k+1}(g)$ (which means that $\lambda_{k}(g)$ corresponds to the first one or the last one of a cluster of equal eigenvalues), then the necessary conditions above are also sufficient.
\end{itemize}

The condition (i) above means that the map $u:=(u_{1}, \cdots , u_{d}): (M,g)\longrightarrow \R^{d}$, is an isometric immersion, whose image is a minimally immersed submanifold of the Euclidean sphere $\mathbb{S}^{d-1}(\sqrt{n \over \lambda_{k}(g)})$ of radius $ \sqrt{n \over \lambda_{k}(g)}$ (see \cite{T}). In other words, a metric $g$ is critical for the functional $\lambda_{k}$, for some $k \geq 1$, if and only if $g$ is induced on $M$ by a minimal immersion of $M$ into a sphere. Therefore, the classification of critical metrics of eigenvalue functionals on $M$ reduces to the classification of minimal immersions of $M$ into spheres.
The many existence and classification results (see for instance \cite{FP,Hs,ER,To} and the references therein) of minimal immersions into spheres give examples of critical metrics for the eigenvalue functionals.  

Notice that, for the first eigenvalue functional, the critical metrics are classified on surfaces of genus $0$ and $1$. Indeed, on $\mathbb{S}^{2}$ and $\R P^{2}$,
the standard metrics are the only critical ones (\cite{EI0}). On the torus $\mathbb{T}^{2}$, the flat metrics $g_{eq}$ and $g_{sq}$ mentioned above are, up to dilatations, the only critical metrics for $\lambda_{1}$ (\cite {EI1}). The metric $g_{eq}$ corresponds to a maximizer for $\lambda_{1}$ (\cite{N1}), while $g_{sq}$ is a saddle point. For the Klein bottle $\mathbb{K}$, Jacobson, Nadirashvili and Polterovitch \cite{JNP} showed the existence of a critical metric and 
El Soufi, Giacomini and Jazar \cite{EGJ} proved that this metric is, up to dilatations, the unique critical metric for $\lambda_{1}$ on $\mathbb{K}$.

Now, the condition (ii), concerning the criticality for $\lambda_{k}$ restricted to a conformal class, is equivalent to the fact that the map $u:=(u_{1}, \cdots , u_{d}): (M,g)\longrightarrow \mathbb{S}^{d-1}$ is a harmonic map with energy density $e(u)= {\lambda_{k}(g) \over 2}$ (see for instance \cite{EL}). Thus, a metric $g$ is critical for some $\lambda_{k}$ restricted to the conformal class of $g$ if and only if $(M,g)$ admits a harmonic map of constant energy density in a sphere. In particular, the metric of any homogeneous compact Riemannian space is critical for $\lambda_{k}$ restricted to its conformal class (for other examples see \cite{PU, ER, To} and the references therein).

A consequence of the necessary condition (ii) is that, if $g$ is a critical metric of $\lambda_{k}$ restricted to $C(g)$, then the multiplicity of $\lambda_{k}(g)$ is at least $2$ (Corollary \ref{cor 42}). This means that $\lambda_{k}(g)=\lambda_{k-1}(g)$ or $\lambda_{k}(g)=\lambda_{k+1}(g)$. In the case where the metric $g$ is a local maximizer of $\lambda_{k}$ restricted to $C(g)$, we prove that one necessarily has: $\lambda_{k}(g)=\lambda_{k+1}(g)$
(Corollary \ref{cor 43}). For a local minimizer, one has $\lambda_{k}(g)=\lambda_{k-1}(g)$. Together with a recent result of Colbois and the first author \cite{CE}, this result tells us that a Riemannian metric can never maximize two consecutive eigenvalues simultaneously on its conformal class (Corollary \ref{cor44}). In fact, if $g$ maximizes $\lambda_{k}$ on $C(g)$, then
$$\lambda_{k+1}(g)^{n \over 2} \leq \sup_{g' \in C(g)} \lambda_{k+1}(g')^{n \over 2} - n^{n \over 2} \omega_{n},$$
where $\omega_{n}$ is the volume of the unit Euclidean n-sphere.\\

As an application of the results above, one can derive characterizations of the metrics which are critical for various functions of eigenvalues. To illustrate this, we treat in the last section of this paper, the case of the ratio functional ${\lambda_{k+1} \over \lambda_{k}}$ of two consecutive eigenvalues and give characterizations of critical metrics for these functionals.


\section {Derivatives of eigenvalues with respect to metric deformations}
Let $M$ be a compact smooth manifold of dimension $n \geq 2$. For any Riemannian metric $g$ on $M$, we denote by $ 0 < \lambda_{1}(g) \leq \lambda_{2}(g)\leq 
\cdots $  the eigenvalues of the Laplace-Beltrami operator $ \Delta_{g}$ associated with $g$. For any $k \in \N$, we denote by $E_{k}(g)= \hbox{Ker}(\Delta_{g}- \lambda_{k}(g)I)$ the eigenspace corresponding to $\lambda_{k}(g)$ and by $ \Pi_{k}: L^{2}(M,g)\longrightarrow E_{k}(g)$ the orthogonal  projection on $E_{k}(g)$.\\
Let us fix a positive integer $k$ and consider the functional $g\longmapsto \lambda_{k}(g)$. This functional is continuous but not differentiable in general. However, perturbation theory tells us that $\lambda_{k}$ is left and right differentiable along any analytic curve of metrics. The main purpose of this section is to express the derivatives of $\lambda_{k}$ with respect to analytic metric deformations, in terms of the eigenvalues of an explicit quadratic form on $E_k(g)$. Indeed, we have the following

\begin{theorem}\label{main} Let $g$ be a Riemannian metric on $M$ and let $(g_{t})_{t}$ be a family of Riemannian metrics analytically indexed by $t \in (-\epsilon,\epsilon)$, such that $g_{0}=g$. The following hold
\begin{enumerate}
\item [i)] The function $t\in (-\epsilon,\epsilon) \longmapsto \lambda_{k}(g_{t})$ admits a left and a right derivatives at $t=0$.
\item [ii)]The derivatives ${d \over dt}\lambda_{k}(g_{t})\big|_{t=0^-}$ and ${d \over dt}\lambda_{k}(g_{t})\big|_{t=0^+}$ are eigenvalues of the operator
$\Pi_{k} \circ \Delta': E_{k}(g) \longrightarrow E_{k}(g)$, where $\Delta'= {d \over dt} \Delta_{g_{t}} \big|_{t=0}$.
\item [iii)]If $\lambda_{k}(g)>\lambda_{k-1}(g)$, then ${d \over dt}\lambda_{k}(g_{t})\big|_{t=0^-}$ and ${d \over dt}\lambda_{k}(g_{t})\big|_{t=0^+}$ are the greatest and the least eigenvalues of $\Pi_{k} \circ \Delta'$ on $E_{k}(g)$, respectively.
\item [iv)]If $\lambda_{k}(g)<\lambda_{k+1}(g)$, then ${d \over dt}\lambda_{k}(g_{t})\big|_{t=0^-}$ and ${d \over dt}\lambda_{k}(g_{t})\big|_{t=0^+}$ are the least and the greatest eigenvalues of $\Pi_{k} \circ \Delta'$ on $E_{k}(g)$, respectively.
\end{enumerate}
\end{theorem}

\begin{proof}
The family of operators $\Delta_{t}:= \Delta_{g_{t}}$ depends analytically on $t$ and, $\forall t, \, \Delta_{t}$ is self-adjoint with respect to the $L^{2}$ inner product induced by $g_{t}$ (but not necessarily to that induced by $g$). However, as done in  \cite{BU}, after a conjugation by the unitary isomorphism
\begin {align*}
U_{t}: L^{2}(M,g) & \to  L^{2}(M,g_{t})\\
{}        u       & \mapsto { \left(\frac{|g|}{|g_{t}|}\right)^{1/4}u},
\end{align*}
where $|g_{t}|$ is the Riemannian volume density of $g_{t}$, we obtain an analytic family $P_{t}=U_{t}^{-1}\circ \Delta_{t}\circ U_{t}$ of operators such that,
$\forall t\in (-\epsilon,\epsilon)$, $P_{t}$ is self-adjoint with respect to the $L^{2}$ inner product induced by $g$. Moreover, $P_{t}$ and $\Delta_{t}$ have the same spectrum. In particular,
$\lambda_{k}(g_t)$ is an eigenvalue of $P_{t}$. The Rellich-Kato perturbation theory of unbounded self-adjoint operators applies to the analytic family of operators $t \mapsto P_{t}$. Therefore, if we denote by $m$ the dimension of $E_{k}(g)$, then there exist, $\forall t\in (-\epsilon,\epsilon)$, $m$ eigenvalues $\Lambda_{1}(t), \cdots ,\Lambda_{1}(t)$ of $P_{t}$ associated with an $L^{2}(M,g)$-orthonormal family of eigenfunctions $v_{1}(t), \cdots , v_{m}(t)$ of $P_{t}$, that is $P_{t} v_{i}(t)= \Lambda_{i}(t) v_{i}(t)$, so that $\Lambda(0)= \cdots = \Lambda_{m}(0)=\lambda_{k}(g)$, and $\forall i \leq m$, both $\Lambda_{i}(t)$ and $v_{i}(t)$ depend analytically on $t$. Setting, $\forall i \leq m$ and $\forall t \in (-\epsilon,\epsilon)$, $u_{i}(t)=U_{t} v_{i}(t)$, we get, $\forall i \leq m$,
\begin{equation}\label{1}
\Delta_{t} u_{i}(t) = \Lambda_{i}(t) u_{i}(t)
\end{equation}
and the family $\left\{u_{1}(t), \cdots, u_{m}(t)\right\}$ is orthonormal in $L^{2}(M,g_{t})$.
Since $t\mapsto \lambda_{k}(t)$ is continuous and, $\forall i \leq m, \; t\mapsto \Lambda_{i}(t)$ is analytic with $\Lambda_{i}(0)=\lambda_{k}(g)$, there exist $\delta >0$ and two integers $p, \, q \leq m$ such that
\[
\lambda_{k}(g_{t})= \left\{ 
\begin{array}{rclll}
\Lambda_{p}(t) & \hbox {for}\; t\in (-\delta,0)\\
\Lambda_{q}(t) & \hbox {for} \; t\in (0,\delta).  
\end{array}
\right.
\]
Assertion $(i)$ follows immediately.
Moreover, one has 
$$ {d \over dt}\lambda_{k}(t) \big|_{t=0^-}= \Lambda'_{p}(0)$$
and
$${d \over dt}\lambda_{k}(t) \big|_{t=0^+}= \Lambda'_{q}(0).$$
Differentiating both sides of (\ref{1}) at $t=0$, we get 
$$ \Delta' u_{i}+\Delta u'_{i}= \Lambda'_{i}(0) u_{i}+\lambda_{k}(g) u'_{i}$$
with $u'_{i}= {d \over dt}u_{i}(t)\big|_{t=0}$ and $u_{i}:= u_{i}(0)$.
Multiplying this last equation by $u_{j}$ and integrating by parts with respect to the Riemannian volume element $v_g$ of $g$, we obtain

\begin{equation} \label{2}
\int_{M}u_{j} \Delta'u_{i} \ v_{g}= \left\{
\begin{array}{cl}
\Lambda'_{i}(0) & \hbox {if} \; j=i \\
  0             & \hbox{otherwise}.
\end{array} \right.
\end{equation}
Since $\left\{u_{1}, \cdots, u_{m}\right\}$ is an orthonormal basis of $E_{k}(g)$ with respect to the $L^{2}$-inner product induced by $g$, we deduce that
$$ (\Pi_{k} \circ \Delta') u_{i}= \Lambda'_{i}(0) u_{i}.$$
In particular, $\Lambda'_{p}(0)$ and $\Lambda'_{q}(0)$ are eigenvalues of $\Pi_{k} \circ \Delta'$, which proves Assertion (ii).

Assume now $ \lambda_{k}(g)> \lambda_{k-1}(g)$. Hence, $\forall i\le m$, $ \Lambda_{i}(0)= \lambda_{k}(g)> \lambda_{k-1}(g)$. By continuity, we have $\Lambda_{i}(t) >\lambda_{k-1}(g_{t})$ for sufficiently small $t$. Hence, there exists $\eta >0$ such that, $ \forall t \in (-\eta,\eta)$ and $\forall i\le m$, $ \Lambda_{i}(t) \geq \lambda_{k}(g_{t})$, which means that $\lambda_{k}(g_{t})= \min\left\{\Lambda_{1}(t), \cdots, \Lambda_{m}(t)\right\}.$
This implies that
$$ {d \over d t} \lambda_{k}(g_{t})\big|_{t=0^-}= \max \left\{\Lambda'_{1}(0), \cdots, \Lambda'_{m}(0)\right\}$$
and
$$ {d \over d t} \lambda_{k}(g_{t})\big|_{t=0^+}= \min \left\{\Lambda'_{1}(0), \cdots, \Lambda'_{m}(0)\right\}.$$
Assertion (iii) is proved.

The proof of Assertion (iv) is similar. Indeed, if $ \lambda_{k}(g)< \lambda_{k+1}(g)$, one has, for sufficiently small $t$, $\lambda_{k}(g_{t})= \max\left\{\Lambda_{1}(t), \cdots, \Lambda_{m}(t)\right\}$ and, then,
$$ {d \over dt} \lambda_{k}(g_{t})\big|_{t=0^+}= \max \left\{\Lambda'_{1}(0), \cdots, \Lambda'_{m}(0)\right\}$$
and
$$ {d \over dt} \lambda_{k}(g_{t})\big|_{t=0^-}= \min \left\{\Lambda'_{1}(0), \cdots, \Lambda'_{m}(0)\right\}.$$

\end{proof}

The quadratic form associated with the symmetric operator $\Pi_{k} \circ \Delta'$ acting on $E_k(g)$ can be expressed explicitely as follows:

\begin{lemma} \label{lem 21}
Let $(g_{t})_{t}$ be an analytic deformation of the metric $g$ and let $h:= {d \over dt} g_{t}\big|_{t=0}$. The operator
$P_{k,h}:= \Pi_{k} \circ \Delta': E_k(g)\to E_k(g)$ is a symmetric with respect to the $L^2$-norm induced by g ; the corresponding quadratic form is given by, $\forall u \in E_{k}(g)$,
$$  Q_{h}(u):= \int_{M} u \, P_{k,h}u \, v_{g}=-\int_M\left( du\otimes du+ {1 \over 4} \Delta_{g}u^{2}g,h\right) v_g,$$
where $(\,,\,)$ is the pointwise inner product induced by $g$ on covariant 2-tensors.
Moreover, if , $\forall t \in (-\epsilon,\epsilon), \, g_{t}=\alpha_{t}g$ is conformal to $g$, then $h=\varphi g$ with $\varphi= {d \over dt}\alpha_{t}\big|_{t=0}$, and, $\forall u \in E_{k}(g)$,
$$ Q_{h}(u)=-\int_{M} \varphi (|du|^{2}+{n \over 4}\Delta_g u^{2}) v_{g}=-{n \over 2}\int_{M} \varphi (\lambda_{k}(g) u^{2}- {n-2 \over n} |du|^{2}) v_{g}.$$
\end{lemma}
\begin{proof}
The derivative at $t=0$ of $t\mapsto \Delta_{g_t} $ is given by the formula (see \cite{B}):
\begin{equation}
\Delta' u:= \frac{d}{dt} \Delta_{g_t} u \Big|_{t=0}= ( Ddu, h) - ( du, \delta h +\frac{1}{2} d(\hbox{trace}_g h)),
\end{equation}
where $D$ is the canonical covariant derivative induced by $g$.
Thus,
\begin{equation}\label{3}
\int_M u \, \Delta' u\; v_g=\int_M u ( Ddu, h) v_g -{1\over 2} \int_M( du^2, \delta h +\frac{1}{2} d(\hbox{trace}_g h) ) v_g.
\end{equation}
One has, $\forall u$,
$$ uDdu =  {1\over 2}Ddu^2- du\otimes du.$$
Hence,
$$\int_M u \left( Ddu, h\right) v_g = {1\over 2}\int_M ( Ddu^2, h) v_g - \int_M  ( du\otimes du, h) v_g.$$
Since $\delta$ is the adjoint of $D$ w.r.t. the $L^2(g)$-inner product, we obtain
$$\int_M u \left( Ddu, h\right) v_g = {1\over 2}\int_M ( du^2, \delta h) v_g - \int_M  ( du\otimes du, h) v_g.$$
On the other hand
$$\int_M( du^2,  d(\hbox{trace}_g h) ) v_g = \int_M \Delta_g u^2 \; \hbox{trace}_g h v_g =\int_M ( \Delta_g u^2  \; g,  h ) v_g .
$$
Replacing in (4) one immediately gets the desired identity.

A straightforward computation gives the expression of $Q_h$ for conformal deformations.
\end{proof}

In relation to the question (2) of the introduction, we give the following result which is a direct consequence of Theorem \ref{main} and Lemma \ref{lem 21}.
\begin{corollary}\label{cor 21}
Let $(g_{t})_{t}$ be an analytic deformation of a Riemannian metric $g$ on $M$ and let $Q_h$ be the associated quadratic form defined as in Lemma \ref{lem 21}, with $h= {d \over dt} g_{t}\big|_{t=0}$.
\begin{itemize}
	\item [i)] If $Q_h$ is positive definite on $E_k(g)$, then there exists $\varepsilon >0$ such that  $\lambda_{k}(g_{-t})<\lambda_{k}(g)<\lambda_{k}(g_{t})$ for all $t\in(0,\varepsilon)$.
	\item [ii)] Assume that $\lambda_{k}(g)>\lambda_{k-1}(g)$. If there exists $u\in E_k(g)$ such that $Q_h(u)<0$, then $\lambda_{k}(g_{t})<\lambda_{k}(g)$ for all $t\in(0,\varepsilon)$, for some  $\varepsilon >0$.
	\item [iii)] Assume that $\lambda_{k}(g)<\lambda_{k+1}(g)$. If there exists $u\in E_k(g)$ such that $Q_h(u)>0$, then $\lambda_{k}(g_{t})>\lambda_{k}(g)$ for all $t\in(0,\varepsilon)$, for some  $\varepsilon >0$.
\end{itemize}

In particular, if $Q_h(u)>0$ for a first eigenfunction $u$, then $\lambda_{1}(g_{t})<\lambda_{1}(g)$ for sufficiently small positive $t$.
\end{corollary}

\section{Critical metrics of the eigenvalue functionals}

Let $M$ be a closed manifold of dimension $n \geq 2$ and let $k$ be a positive integer. Before introducing the notion of critical metric of the functional $\lambda_{k}$, notice that this functional is not scaling invariant. Therefore, we will restrict $\lambda_{k}$ to the set of metrics of given volume.
In view of Theorem \ref{main}, a natural way to introduce the notion of critical metric is the following:
\begin{definition}
A metric $g$ on $M$ is said to be ``critical'' for the functional $\lambda_{k}$ if, for any volume-preserving analytic deformation $(g_{t})_{t}$ of $g$ with
$g_{0}=g$, the left and the right derivatives of $\lambda_{k}(g_{t})$ at $t=0$ satisfy
$$ {d \over dt}\lambda_{k}(g_{t})\big|_{t=0^-} \times {d \over dt}\lambda_{k}(g_{t})\big|_{t=0^+} \leq 0. $$
\end{definition}

It is easy to see that 
$$ {d \over dt}\lambda_{k}(g_{t})\big|_{t=0^+} \leq 0 \leq {d \over dt}\lambda_{k}(g_{t})\big|_{t=0^-}\Longleftrightarrow \lambda_{k}(t) \leq \lambda_{k}(0)+o(t)$$
and
$$ {d \over dt}\lambda_{k}(g_{t})\big|_{t=0^-} \leq 0 \leq {d \over dt}\lambda_{k}(g_{t})\big|_{t=0^+}\Longleftrightarrow \lambda_{k}(t) \geq \lambda_{k}(0)+o(t).$$
Therefore, $g$ is critical for $\lambda_{k}$ if, for any volume-preserving analytic deformation $(g_{t})_{t}$ of $g$, one of the following inequalities holds:
$$ \lambda_{k}(g_{t}) \leq \lambda_{k}(g) +o(t)$$
or
$$ \lambda_{k}(g_{t}) \geq \lambda_{k}(g) +o(t).$$
Of course, if $g$ is a local maximizer or a local minimizer of $\lambda_{k}$, then $g$ is critical in the sense of the previous definition.

In all the sequel, we will denote by $S^{2}_{0}(M,g)$ the space of covariant $2$-tensors $h$ satisfying $\int_{M} \hbox {trace}_{g}h v_{g}= \int_{M} (g,h) v_{g}=0$, endowed with its natural $L^2$ norm induced by $g$.
\begin{proposition}\label{prop 31}
If $g$ is a critical metric for the functional $\lambda_{k}$ on $M$, then, $\forall h \in S^{2}_{0}(M,g)$, the quadratic form 
$$  Q_{h}(u)=-\int_M\left( du\otimes du+ {1 \over 4} \Delta_{g}u^{2}g,h\right) v_g$$
 is indefinite on $E_{k}(g).$
\end{proposition}
\begin{proof}
Let $h \in S^{2}_{0}(M,g)$. The deformation of $g$ defined for small $t$ by $ g_{t}=\displaystyle{\left[{\mbox{vol}(g) \over \mbox{vol}(g+th)}\right]^{2/n}(g+th)}$, where $\mbox{vol}(g)$ is the Riemannian volume of $(M,g)$, is volume-preserving and depends analytically on $t$ with ${d \over dt}g_{t}\big|_{t=0}=h$. Using Theorem \ref{main}, we see that, if $g$ is critical, then the operator $P_{h,k}$ admits a nonnegative and a nonpositive eigenvalues on $E_{k}(g)$ which means that the quadratic form $Q_{h}$ is indefinite (Lemma \ref{lem 21}).
\end{proof}
In the case where $\lambda_{k}(g)>\lambda_{k-1}(g)$ or $\lambda_{k}(g)<\lambda_{k+1}(g)$, one can show that the converse of Proposition \ref{prop 31} is also true. Indeed, we have the following
\begin{proposition}\label{prop 32}
Let $g$ be a Riemannian metric on $M$ such that $\lambda_{k}(g)>\lambda_{k-1}(g)$ or $\lambda_{k}(g)<\lambda_{k+1}(g)$.
Then $g$ is critical for the functional $\lambda_{k}$ if and only if, $\forall h \in S^{2}_{0}(M,g)$, the quadratic form $Q_{h}$ is indefinite on $E_{k}(g)$.
\end{proposition}
\begin{proof}
Let $(g_{t})_{t}$ be an analytic volume-preserving deformation of $g$ and let $h={d \over dt} g_{t}\big|_{t=0}$. Since $\mbox{vol}(g_{t})$ is constant with respect to $t$, the tensor $h$ belongs to $S^{2}_{0}(M,g)$ (indeed, $\int_{M} (g,h) v_{g}={d \over dt} \hbox{vol}(g_{t})\big|_{t=0}=0$). The indefiniteness of $Q_{h}$ implies that the operator $P_{k,h}=\Pi_{k} \circ \Delta'$ admits both non-negative and non-positive eigenvalues on $E_{k}(g)$ (see Lemma \ref{lem 21}). The result follows immediately from Theorem \ref{main} (iii) and (iv).
\end{proof}
The indefiniteness of $Q_{h}$ on $E_{k}(g)$ for all $h \in S^{2}_{0}(M,g)$, can be interpreted intrinsically in terms of the eigenfunctions of $\lambda_{k}(g)$ as follows.
\begin{lemma}\label {lem 31}
Let $g$ be a Riemannian metric on $M$. The two following conditions are equivalent:
\begin{enumerate}
\item[i)] $\forall h \in S^{2}_{0}(M,g)$, the quadratic form $Q_{h}$ is indefinite on $E_{k}(g)$.
\item[ii)] There exists a finite family $\left\{u_{1}, \cdots , u_{d}\right\}\subset E_{k}(g)$ of eigenfunctions associated with $\lambda_{k}(g)$ such that
$$ \sum_{i \leq d} du_{i} \otimes du_{i}=g.$$
\end{enumerate}
\end{lemma} 

\begin{proof}
The proof of ``(i) implies (ii)'' uses the same arguments as in the proof of Theorem 1.1 of  \cite{EI1}. For the sake of completeness, we will recall the main steps. First, we introduce the convex set $K \subset S^{2}(M,g)$ given by 
$$K= \left\{\sum_{j \in J} \left[du_{j}\otimes du_{j} + {1 \over 4} \Delta_{g}u_{j}^{2}\, g\right]; u_{j}\in E_{k}(g), \, J \subset \N, \, J \,\hbox{finite} \right\}.$$
Let us first show that $g\in K$. Indeed, if $g \notin K$, then, applying classical separation theorem in the finite dimensional subspace of $S^{2}(M,g)$ generated by $K$ and $g$, endowed with the $L^{2}$ inner product induced by $g$, we deduce the existence of a $2$-tensor $h\in S^{2}(M,g)$ such that $\int_{M}(g,h) v_{g}>0$ and, $\forall\ T \in K, \int_{M} (T,h) v_{g} \leq 0. $
The tensor  $$h_{0}= h- \left(\displaystyle{{1 \over n \ \mbox{vol}(g)}}\int_{M} (g,h) v_{g}\right)g$$  belongs to $S^{2}_{0}(M,g)$ and we have, $\forall u \in E_{k}(g)$, $u \ne 0$,
$$Q_{h_{0}}(u)= - \int_{M} (du \otimes du+{1 \over 4}\Delta_{g}u^{2}g,h)v_{g} + {\int_{M}(g,h)v_{g} \over n \,\mbox{vol}(g)} \int_{M} |du|^{2} v_{g}$$
$$\geq {\lambda_{k}(g) \over n \,\mbox{vol}(g)} \int_{M} (g,h)v_{g} \int_{M}u^{2} v_{g}.$$
Since $\int_{M} (g,h) v_{g} > 0$, the quadratic form $Q_{h_{0}}$ is positive definite, which contradicts the assumption (i). Now, $g\in K$ means that there exists $u_{1}, \cdots , u_{m} \in E_{k}(g)$ such that
\begin{equation}\label{5}
 \sum_{i \leq d} (du_{i} \otimes du_{i} + {1 \over 4} \Delta_{g}u_{i}^{2})g=g.
 \end{equation}
Hence, since $\Delta u_{i}^{2}=2(\lambda_{k}(g)u_{i}^{2}-|du_{i}|^{2})$, we obtain after taking the trace in (\ref{5}),
$${\lambda_{k}(g) \over 2}\sum_{i \leq d}u_{i}^{2}= 1+ {n-2 \over 2n} \sum_{i \leq d} |du_{i}|^{2}.$$
For $n=2$, we immediately get $\sum_{i \leq d} u_{i}^{2}= {2 \over \lambda_{k}(g)}$ and, for $n\ge 3$, we consider the function $f:=\sum_{i \leq d}u_{i}^{2} - {n \over \lambda_{k}(g)}$ and observe that it satisfies 
$$ (n-2) \Delta_{g}f= 2(n-2)(\lambda_{k}(g) \sum_{i \leq d} u_{i}^{2} - \sum_{i \leq d}|du_{i}|^{2})=-4\lambda_{k}(g) f.$$
Thus, $f=0$ (the Laplacian being a non-negative operator) and, then, $\forall n\ge 2$, $\sum_{i \leq d} u_{i}^{2}= {n \over \lambda_{k}(g)}$. Replacing in (\ref{5}), we obtain 
$$ \sum_{i \leq d} du_{i} \otimes du_{i}=g.$$
Conversely, let $u_{1}, \cdots, u_{d}$ be as in (ii). This means that the map $x\in M \longmapsto u(x)=(u_{1}(x), \cdots , u_{d}(x)) \in \R^{d}$ is an isometric immersion. The vector $\Delta u(x)= (\Delta u_{1}(x),\cdots , \Delta u_{d}(x))= \lambda_{k}(g) u(x)$ represents the mean curvature vectorfield of the immersed submanifold $u(M)$. Hence, $\forall x \in M$, the position vector $u(x)$ is normal to $u(M)$ which implies that $u(M)$ is contained in a sphere of $\R^{d}$ centered at the origin. Thus, $\sum_{i \leq d}u_{i}^{2}$ is constant on $M$. Consequently, $\forall h \in S^{2}_{0}(M,g),$
$$ \sum_{i \leq d} Q_{h}(u_{i})= \cdots =0.$$
It follows that $Q_{h}$ is indefinite on $E_{k}(g)$. 
\end{proof}

Propositions \ref{prop 31}, \ref{prop 32} and Lemma \ref {lem 31} lead to the following characterization of critical metrics of $\lambda_{k}$:
\begin{theorem} \label {th 31}
Let $g$ be a Riemannian metric on $M$.
\begin{enumerate}
\item[i)] If $g$ is a critical metric of the functional $\lambda_{k}$, then there exists a finite family $\left\{u_{1}, \cdots , u_{d}\right\} \subset E_{k}(g)$
of eigenfunctions associated with $\lambda_{k}(g)$ such that
$$ \sum_{i \leq d} du_{i} \otimes du_{i} =g. $$
\item[ii)] Assume that $\lambda_{k}(g)> \lambda_{k-1}(g)$ or $\lambda_{k}(g)< \lambda_{k+1}(g)$. Then $g$ is a critical metric of the functional $\lambda_{k}$
if and only if there exists a finite family $\left\{u_{1}, \cdots , u_{d}\right\} \subset E_{k}(g)$ of eigenfunctions associated with $\lambda_{k}(g)$ such that
$$ \sum_{i \leq d} du_{i} \otimes du_{i} =g. $$
\end{enumerate}
\end{theorem}
According to Theorem \ref{th 31} (ii), the standard metrics $g$ of compact rank one symmetric spaces are critical metrics of the functionals $\lambda_{k}$, for any $k$ such that 
$\lambda_{k}(g)> \lambda_{k-1}(g)$ or $\lambda_{k}(g)< \lambda_{k+1}(g)$. More generally, this is the case of all compact Riemannian homogeneous spaces with irreducible isotropy representation. Indeed, if $\left\{u_{1}, \cdots , u_{d}\right\}$ is an $L^{2}(g)$-orthonormal basis of $E_{k}(g)$, then the tensor $\sum_{i \leq d}du_{i} \otimes du_{i}$ is invariant under the isometry group action which implies that it is proportional to $g$ (Schur's Lemma).

In \cite{EI3}, we studied the notion of critical metrics of the trace of the heat kernel $Z_{g}(t)= \sum e^{-\lambda_{k}(g) t}$, considered as a functional on the set of metrics of given volume. We obtain various characterizations of these critical metrics. An immediate consequence of Theorem \ref{th 31} and  \cite[Theorem 2.2]{EI3},  is the following
\begin{corollary}
Let $g$ be a Riemannian metric on $M$. If $g$ is a critical metric of the trace of the heat kernel at any time $t>0$, then $g$ is a critical metric of the functional $\lambda_{k}$ for all $k$ such that $\lambda_{k}(g)> \lambda_{k-1}(g)$ or $\lambda_{k}(g)< \lambda_{k+1}(g)$. 
\end{corollary}
In particular, the flat metrics $g_{sq}$ and $g_{eq}$ on the 2-Torus $\mathbb{T}^{2}$ associated with the square lattice $\Z^{2}$ and the equilateral lattice
$\Z (1,0)\oplus \Z (1/2,\sqrt{3}/2)$, respectively, are critical metrics of the functionals $\lambda_{k}$ for all $k$ such that $\lambda_{k}(g_{sq})> \lambda_{k-1}(g_{sq})$ or $\lambda_{k}(g_{sq})< \lambda_{k+1}(g_{sq})$ and such that $\lambda_{k}(g_{eq})> \lambda_{k-1}(g_{eq})$ or $\lambda_{k}(g_{eq})< \lambda_{k+1}(g_{eq})$ respectively (see \cite{EI3}). Other examples of critical metrics can be obtained as Riemannian products of previous examples (see \cite{EI3}).

As we noticed in the proof of Lemma \ref {lem 31}, the condition $ \sum_{i \leq d} du_{i} \otimes du_{i} =g$, with $u_{i} \in E_{k}(g)$, implies that the map $u=(u_{1}, \cdots , u_{d})$ is an isometric immersion of $(M,g)$ into a $(d-1)$-dimensional sphere. In particular, the rank of $u$ is at least $n$. Therefore we have the following
\begin{corollary}
If $g$ is a critical metric of the functional $\lambda_{k}$, then 
$$ \dim E_{k}(g) \geq \dim M+1.$$
Moreover, the equality implies that $(M,g)$ is isometric to an Euclidean sphere. 
\end{corollary}
In the particular case where a metric $g$ is a local maximizer of $\lambda_{k}$ (that is $\lambda_{k}(g_{t})\leq \lambda_{k}(g)$ for any volume-preserving deformation $(g_{t})_{t}$ of $g$), we have the additional necessary condition that $ \lambda_{k}(g)=\lambda_{k+1}(g)$ (see Proposition \ref{cor 43} below). For a local minimizer, we have 
$ \lambda_{k}(g)=\lambda_{k-1}(g)$. In particular, the functional $\lambda_{1}$ admits no local minimizing metric. This result have been obtained by us in \cite{EI1} using different arguments.

We end this section with the following result in the spirit of Berger's work \cite{B}:
\begin{corollary} \label{cor33}
Let $g$ be a Riemannian metric on $M$. Let $p \geq 1$ and $q \geq p$ be two natural integers such that 
$$ \lambda_{p-1}(g) < \lambda_{p}(g)=\lambda_{p+1}(g)= \dots = \lambda_{q}(g)< \lambda_{q+1}(g).$$
The metric $g$ is critical for the functional $ \sum_{i=p}^{q}\lambda_{i}$ if and only if there exists an $L^{2}(M,g)$-orthonormal basis $u_{1}, u_{2}, \dots , u_{m}$ of $E_{p}(g)$ such that $\sum_{i=1}^{m}du_{i} \otimes du_{i}$ is proportional to $g$.
\end{corollary}
\begin{proof}
The multiplicity of $\lambda_{p}(g)$ is $m=q-p+1$. Let $(g_{t})_{t}$ be a volume-preserving analytic deformation of $g$ and $h={d \over dt} g_{t}\big|_{t=0} \in S_{0}^{2}(M,g)$. Let $\Lambda_{1}(t), \cdots, \Lambda_{m}(t)$ and $v_{1}(t), \cdots, v_{m}(t)$ be the families of eigenvalues and orthonormal eigenfunctions of $\Delta_{g_{t}}$ depending analytically on $t$ and such that $\Lambda_{1}(0)= \cdots = \Lambda_{m}(0)=\lambda_{p}(g)$, as in the proof of Theorem \ref{main}. For sufficiently small $t$, one has
$$ \sum_{i=p}^{q}\lambda_{i}(g_{t})=\sum_{i=1}^{m}\Lambda_{i}(t).$$
Hence, $\sum_{i=p}^{q}\lambda_{i}(g_{t})$ is differentiable at $t=0$ and one has (see the proof of Theorem \ref{main} and Lemma \ref{lem 21}),
$ \frac{d}{dt} \sum_{i=p}^{q}\lambda_{i}(g_{t})\big|_{t=0} = \sum_{i=1}^{m}\Lambda'_{i}(0)= \sum_{i=1}^{m} Q_{h}(v_{i})$, with $v_{i}:=v_{i}(0)$. Therefore, $g$ is critical for $\sum_{i=p}^{q}\lambda_{i}$ if and only if, $ \forall h \in S_{0}^{2}(M,g), \sum_{i=1}^{m} Q_{h}(v_{i})=0$. As in the proof of Lemma 3.1, this last condition means that $\sum_{i \le m} dv_{i} \otimes dv_{i}$ is proportional to $g$.
\end{proof}


\section{Critical metrics of the eigenvalue functionals in a conformal class}

Let $M$ be a closed manifold of dimension $n \geq 2$. For any Riemannian metric $g$ on $M$, we will denote by $C(g)$ the set of metrics which are conformal to $g$ and have the same volume as $g$, i.e
$$ C(g)= \left\{e^{\alpha}g \,;\, \alpha \in \mathcal{C}^{\infty}(M)\, \hbox {and} \; vol(e^{\alpha}g)= vol(g)\right\}.$$
Let $k$ be a positive integer. The purpose of this section is to study critical metrics of the functional $\lambda_{k}$
restricted to a conformal class $C(g)$.
\begin{definition}
A metric $g$ is said to be critical for the functional $\lambda_{k}$ restricted to $C(g)$ if, for any analytic deformation $\left\{g_{t}= e^{\alpha_{t}}g\right\}\subset C(g)$ with $g_{0}=g$, we have 
$$ {d \over dt}\lambda_{k}(g_{t})\big|_{t=0^-} \times {d \over dt}\lambda_{k}(g_{t})\big|_{t=0^+} \leq 0.$$
\end{definition}
In the sequel, we denote by $\mathcal{A}_{0}(M,g)$ the set of regular functions $\varphi$ with zero mean on $M$, that is, $\int_{M} \varphi \, v_{g}=0$.
In the spirit of Propositions \ref{prop 31} and \ref{prop 32}, we obtain in the conformal setting, the following
\begin{proposition}\label{prop 41}
Let $g$ be a Riemannian metric on $M$.
\begin{enumerate}
\item[i)] If $g$ is a critical metric of the functional $\lambda_{k}$ restricted to $C(g)$, then, $\forall \varphi \in \mathcal{A}_{0}(M,g)$, the quadratic form
$$ q_{\varphi}(u)= \int_{M} (\lambda_{k}(g)u^{2}-{n-2 \over n}|du|^{2})\varphi\, v_{g}$$
is indefinite on $E_{k}(g)$.
\item[ii)] Assume that $\lambda_{k}(g)> \lambda_{k-1}(g)$ or $\lambda_
{k}(g)< \lambda_{k+1}(g)$. The metric $g$ is critical for the functional $\lambda_{k}$ restricted to $C(g)$ if and only if, $\forall \varphi \in \mathcal{A}_{0}(M,g)$, the quadratic form $q_{\varphi}$ is indefinite on $E_{k}(g)$.
\end{enumerate}
\end{proposition}
\begin{proof}

i) Let $\varphi \in \mathcal{A}_{0}(M,g)$. The conformal deformation of $g$ given by 
	$$ g_{t}:= \left[ {\mbox{vol}(g) \over \mbox{vol}(e^{t\varphi}g)}\right]^{2 \over n} e^{t \varphi}g,$$
belongs to $C(g)$ and depends analytically on $ t $ with $ {d \over dt} g_{t}\big|_{t=0}= \varphi g $.
Following the arguments of the proof of Proposition \ref{prop 31}, we show that the criticality of $g$ for $\lambda_{k}$ restricted to $C(g)$ implies the indefiniteness of  the quadratic form
$Q_{\varphi g}$ on $E_{k}(g)$. Applying Lemma \ref{lem 21}, we observe that $Q_{\varphi g}= -\frac{n}{2}q_{\varphi}$.

ii) Let $g_{t}= e^{\alpha_{t} g} \in C(g)$ be an analytic deformation of $g$. Since $\mbox{vol}(g_{t})$ is constant with respect to $t$, the function $\varphi = {d \over dt}\alpha_{t}\big|_{t=0}$ belongs  $\mathcal{A}_{0}(M,g)$. Applying Theorem \ref{main} (iii) and (iv) and Lemma \ref{lem 21} with $h=\varphi g$, we get the result.
\end{proof}
\begin{lemma}\label{lem 41}
Let $g$ be a Riemannian metric on $M$. The two following conditions are equivalent:
\begin{enumerate}
	\item[i)] $\forall \ \varphi \in \mathcal{A}_{0}(M,g)$, the quadratic form $q_{\varphi}$ is indefinite on $E_{k}(g)$.
	\item[ii)] There exists a finite family $\left\{u_{1}, \cdots , u_{d}\right\} \subset E_{k}(g)$ of eigenfunctions associated with $\lambda_{k}(g)$ such that
	$$ \sum_{i \leq d} u_{i}^{2}=1.$$
\end{enumerate}
\end{lemma}
\begin{proof}
``(i) implies (ii)'': We introduce the convex set 
$$ H=\left\{\sum_{i \in I}\left[\lambda_{k}(g)u_{i}^{2}-{n-2 \over n}|du_{i}|^{2}\right]; \, u_{i} \in E_{k}(g), \, I \subset \N, \, I \,\hbox {finite}\right\}.$$
Using the same arguments as in the proof of Lemma \ref{lem 31}, we show that the constant function $1$ belongs to $H$. Hence, there exist $u_{1}, \cdots , u_{d} \in E_{k}(g)$ such that
$$ \sum_{i \leq d} (\lambda_{k}(g)u_{i}^{2}-{n-2 \over n}|du_{i}|^{2})=\frac{2}{n}\lambda_{k}(g). $$
For $n=2$, we immediately get $\sum_{i \leq d} u_{i}^{2}= 1$. For $n\ge 3$, we set
$f= \sum_{i \leq d}u_i^2 -1 $ and get, after a straightforward calculation,
$$ {n-2 \over 4} \Delta_{g}f= - \lambda_{k}(g) f.$$
Thus, $f=0$ and $\sum_{i \leq d}u_{i}^{2}=1$.

\noindent ``(ii) implies (i)'': let $u_{1}, \cdots , u_{d} \in E_{k}(g)$ be such that $\sum_{i \leq d}u_{i}^{2}=1$. One has
$$ \sum_{i \leq d} |du_{i}|^{2}= -{1 \over 2}\Delta_{g}\sum_{i \leq d}u_{i}^{2}+ \lambda_{k}(g) \sum_{i \leq d}u_{i}^{2}= \lambda_{k}(g).$$
Therefore, $\forall \varphi \in \mathcal{A}_{0}(M,g)$,
$$\sum_{i \leq d} q_{\varphi}(u_{i})= \frac{2}{n}\int_{M} \lambda_{k}(g) \varphi v_{g}=0 $$
which implies the indefiniteness of $q_{\varphi}$.
\end{proof}
Proposition \ref{prop 41} and Lemma \ref{lem 41} lead to the following
\begin{theorem}\label{th 41}
Let $g$ be a Riemannian metric on $M$.
\begin{enumerate}
	\item [i)]If $g$ is a critical metric of the functional $\lambda_{k}$ restricted to $C(g)$, then there exists a finite family $\left\{u_{1}, \cdots , u_{d}\right\} \subset E_{k}(g)$ of eigenfunctions associated with $\lambda_{k}$ such that $\sum_{i \leq d}u_{i}^{2}=1$.
	\item[ii)]Assume that $\lambda_{k}(g)> \lambda_{k-1}(g)$ or $\lambda_{k}(g)< \lambda_{k+1}(g)$. Then, $g$ is critical for the functional $\lambda_{k}$ restricted to $C(g)$ if and only if, there exists a finite family $\left\{u_{1}, \cdots , u_{d}\right\} \subset E_{k}(g)$ of eigenfunctions associated with $\lambda_{k}(g)$ such that
	$$ \sum_{i \leq d} u_{i}^{2}=1.$$
\end{enumerate}
\end{theorem}
The Riemannian metric $g$ of any homogeneous Riemannian space $(M,g)$ is a critical metric of the functional $\lambda_{k}$ restricted to $C(g)$ for all $k$ such that $\lambda_{k}(g)> \lambda_{k-1}(g)$ or $\lambda_{k}(g)< \lambda_{k+1}(g)$. Indeed, any $L^{2}(g)$-orthonormal basis $\left\{u_{i}\right\}_{i \leq d}$ of $E_{k}(g)$ is such that $\sum_{i \leq d}u_{i}^{2}$ is constant on $M$.
In \cite[Theorem 4.1]{EI3}, we proved that a metric $g$ is critical for the trace of the heat kernel restricted to $C(g)$ if and only if its heat kernel $K$ is constant on the diagonal of $M \times M$. This last condition implies that, $\forall k \in \N^{\ast}$, any $L^{2}(g)$-orthonormal basis $\left\{u_{i}\right\}_{i \leq d}$ of $E_{k}(g)$ is such that $\sum_{i \leq d}u_{i}^{2}$ is constant on $M$. Hence, we have the following
\begin{corollary}\label{cor 41}
Let $g$ be a Riemannian metric on $M$ and let $K$ be the heat kernel of $(M,g)$. Assume that, $\forall t>0$, the function $x \in M \longmapsto K(t,x,x)$ is constant, then the metric $g$ is critical for the functional $\lambda_{k}$ restricted to $C(g)$ for all $k$ such that $\lambda_{k}(g)> \lambda_{k-1}(g)$ or $\lambda_{k}(g)< \lambda_{k+1}(g)$.
\end{corollary}
An immediate consequence of Theorem \ref{th 41} is the following
\begin{corollary}\label{cor 42}
If $g$ is a critical metric of the functional $\lambda_{k}$ restricted to $C(g)$, then $\lambda_{k}(g)$ is a degenerate eigenvalue, that is 
$$ \dim E_{k}(g) \geq 2.$$
\end{corollary}
This last condition means that at least one of the following holds: $\lambda_{k}(g)= \lambda_{k-1}(g)$ or $\lambda_{k}(g)= \lambda_{k+1}(g)$. In the case when $g$ is a local maximizer or a local minimizer, we have the following more precise result
\begin{proposition}\label{cor 43}
\begin{enumerate}
\item[i)]If $g$ is a local maximizer of the functional $\lambda_{k}$ restricted to $C(g)$, then $\lambda_{k}(g)=\lambda_{k+1}(g)$.  
\item[ii)]If $g$ is a local minimizer of the functional $\lambda_{k}$ restricted to $C(g)$, then $\lambda_{k}(g)=\lambda_{k-1}(g)$.
\end{enumerate}
\end{proposition}
\begin{proof}
Assume that $g$ is a local maximizer and that $\lambda_{k}(g)< \lambda_{k+1}(g)$. Let $\varphi \in \mathcal{A}_{0}(M,g)$ and let $g_{t}=e^{\alpha_{t}}g \in C(g)$ be a volume-preserving analytic deformation of $g$ such that ${d \over dt}g_{t}\big|_{t=0}= \varphi g$. Denote by $\Lambda_{1}(t), \cdots, \Lambda_{m}(t)$,
with $m=\dim E_{k}(g)$, the associated family of eigenvalues of $\Delta_{g_{t}}$, depending analytically on $t$ and such that $\Lambda_{1}(0)= \cdots= \Lambda_{m}(0)=\lambda_{k}(g)$ (see the proof of Theorem \ref{main}). For continuity reasons, we have, for sufficiently small $t$ and all $i \leq m$,
$$\Lambda_{i}(t) < \lambda_{k+1}(g_{t}).$$
Hence, $\forall i \leq m$ and $\forall t$ sufficiently small,
$$ \Lambda_{i}(t) \leq \lambda_{k}(t) \leq \lambda_{k}(g)=\Lambda_{i}(0).$$
Consequently, $\Lambda'_{i}(0)=0$ for all $i \leq m$. Since $\Lambda'_{1}(0), \cdots , \Lambda'_{m}(0)$ are eigenvalues of the operator $\Pi_{k} \circ \Delta'$
(by Theorem \ref{main}), this operator is identically zero on $E_{k}(g)$. Applying Lemma \ref{lem 21}, we deduce that, $\forall \varphi \in \mathcal{A}_{0}(M,g)$, $
Q_{\varphi g}\equiv 0$ on $E_{k}(g)$. Thus, $\forall u \in E_{k}(g)$, there exists a constant $c$ so that 
$$ |du|^{2}+ {n \over 4} \Delta_{g}u^{2}= c.$$
 Integrating, we get $ c= {\lambda_{k}(g) \over \mbox{vol}(g)} \int_{M} u^{2} v_{g}$. Since $\Delta_{g}u^{2}=2(\lambda_{k}u^{2}-|du|^{2})$,
we obtain
$${n \over 2} u^{2} - {n-2 \over 2\lambda_{k}(g)}|du|^{2}= {1 \over \mbox{vol}(g)}\int_{M} u^{2} v_{g}.$$
Let $x_{0} \in M$ be a point where $u^{2}$ achieves its maximum. At $x_{0}$, we have $du(x_{0})=0$ and, then,
$$ {n \over 2} \max u^{2}= {n \over 2} u^{2}(x_{0})= {1 \over \mbox{vol}(g)}\int_{M} u^{2} v_{g}$$
which leads to a contradiction (since $u$ is not constant and ${n \over 2} \geq 1$).\\
A similar proof works for (ii).  
\end{proof}
In \cite{CE}, Colbois and the first author proved that
\begin{equation}\label{CE}
\sup_{g' \in C(g)} \lambda_{k+1}(g')^{n \over 2} - \sup_{g' \in C(g)} \lambda_{k}(g')^{n \over 2}\geq n^{n \over 2} \omega_{n},
\end{equation}
where $\omega_n$ is the volume of the unit euclidean sphere of dimension $n$. 

An immediate consequence of this result and Proposition \ref{cor 43} is the following
\begin{corollary}\label{cor44}
Let $g$ be a Riemannian metric on $M$. Assume that $g$ maximizes the functional $\lambda_{k}$ restricted to $C(g)$, that is  
$$\lambda_{k}(g)=\sup_{g' \in C(g)} \lambda_{k}(g').$$ 
Then $g$ cannot maximize neither $\lambda_{k+1}$ nor $\lambda_{k-1}$ (for $k \geq 2$) on $C(g)$.
\end{corollary}
More precisely, if $g$ maximizes $\lambda_{k}$ on $C(g)$, then, using  Proposition \ref{cor 43} and (\ref{CE}),
$$\lambda_{k+1}(g)^{n \over 2} \leq \sup_{g' \in C(g)} \lambda_{k+1}(g')^{n \over 2} - n^{n \over 2} \omega_{n}.$$

Finally we have the following conformal version of Corollary \ref{cor33}
\begin{corollary} \label{cor45}
Let $g$ be a Riemannian metric on $M$. Let $p \geq 1$ and $q \geq p$ be two natural integers such that 
$$ \lambda_{p-1}(g) < \lambda_{p}(g)=\lambda_{p+1}(g)= \dots = \lambda_{q}(g)< \lambda_{q+1}(g).$$
The metric $g$ is critical for the functional $ \sum_{i=p}^{q}\lambda_{i}$ restricted to $C(g)$ if and only if there exists an $L^{2}(M,g)$-orthonormal basis $u_{1}, u_{2}, \dots , u_{m}$ of $E_{p}(g)$ such that $\sum_{i=1}^{m}u_{i}^{2}$ is constant on $M$.
\end{corollary}

\section{Critical metrics of the eigenvalue ratios functionals}
Let $M$ be a closed manifold of dimension $n \geq 2$ and let $k$ be a positive integer. This section deals with the functional $ g \longmapsto \displaystyle{{\lambda_{k+1}(g)\over \lambda_{k}(g)}}$. This functional is invariant under scaling, so it is not necessary to fix the volume of metrics under consideration. If $(g_{t})_{t}$ is any analytic deformation of a metric $g$, then $t \longmapsto \displaystyle{{\lambda_{k+1}(g_{t})\over \lambda_{k}(g_{t})}}$ admits left and right derivatives at $t=0$ (Theorem \ref{main}).
\begin{definition}
i) A metric $g$ is said to be critical for the ratio $\displaystyle{{\lambda_{k+1}\over \lambda_{k}}}$ if, for any analytic deformation $(g_{t})$ of $g$, the left and right derivatives of $\displaystyle{{\lambda_{k+1}(g_{t})\over \lambda_{k}(g_{t})}}$ at $t=0$ have opposite signs.\\
ii) The metric $g$ is said to be critical for the ratio functional $\displaystyle{{\lambda_{k+1}\over \lambda_{k}}}$ restricted to the conformal class $C(g)$ if the condition above holds for any conformal analytic deformation $g_{t}=e^{\alpha_{t}}g$ of $g$.
\end{definition}
Let $g$ be a Riemannian metric on $M$. For any covariant 2-tensor $h \in S^{2}(M)$, we introduce the following operator
$$ \tilde{P}_{k,h}: E_{k}(g) \otimes E_{k+1}(g)\longrightarrow E_{k}(g) \otimes E_{k+1}(g)$$
defined by $$ \tilde{P}_{k,h}= \lambda_{k+1}(g) P_{h,k}\otimes Id_{E_{k+1}(g)}-\lambda_{k}(g) Id_{E_{k}(g)}\otimes P_{k+1,h},$$ where $P_{k,h}$ is defined in Lemma \ref{lem 21}. The quadratic form naturally associated with $\tilde{P}_{k,h}$ is denoted by $\tilde{Q}_{k,h}$ and is given by, $\forall u \in E_{k}(g) \, \hbox {and}\, \forall v \in E_{k+1}(g)$,
$$ \tilde{Q}_{k,h}(u\otimes v)=\lambda_{k+1}(g) \left\|v\right\|^{2}_{L^{2}(g)} Q_{h}(u)- \lambda_{k}(g) \left\|u\right\|^{2}_{L^{2}(g)} Q_{h}(v),$$
where $ Q_{h}(u)= -\int_{M} (du \otimes du + {1 \over 4} \Delta_{g}u^{2}g,h) v_{g}.$

Of course, if $\lambda_{k+1}(g)=\lambda_{k}(g)$, then $g$ is a global minimizer of the ratio ${\lambda_{k+1}\over \lambda_{k}}$. Notice that, thanks to Colin de Verdi\`ere's result \cite{CV}, $\forall k\ge 1$, any closed manifold $M$ carries a metric $g$ such that $\lambda_{k+1}(g)=\lambda_{k}(g)$.  A general characterization of critical metrics of ${\lambda_{k+1}\over \lambda_{k}}$ is given in what follows
\begin{proposition}\label{prop 51}
A Riemannian metric $g$ on $M$ is critical for the functional ${\lambda_{k+1}\over \lambda_{k}}$ if and only if, $\forall h \in S^{2}(M)$, the quadratic form $\tilde{Q}_{k,h}$ is indefinite on $ E_{k}(g) \otimes E_{k+1}(g)$.
\end{proposition}
\begin{proof}
The case where $\lambda_{k+1}(g)=\lambda_{k}(g)$ is obvious ($\tilde{Q}_{k,h} (u\otimes u)=0$).  Assume that $\lambda_{k+1}(g)>\lambda_{k}(g)$ and let $(g_{t})_{t}$ be an analytic deformation of $g$. From Theorem \ref{main}, ${d \over dt}\lambda_{k}(g_{t})\big|_{t=0^-}$
 and ${d \over dt}\lambda_{k}(g_{t})\big|_{t=0^+}$ are the least and the greatest eigenvalues of $P_{k,h}$ on $E_{k}(g)$ respectively.\\
Similarly, ${d \over dt}\lambda_{k}(g_{t})\big|_{t=0^-}$ and ${d \over dt}\lambda_{k}(g_{t})\big|_{t=0^+}$ are the greatest and the least eigenvalues of $P_{k+1}$ on $E_{k}(g)$. Therefore,
$$\lambda_{k}(g)^{2} {d \over dt}{\lambda_{k+1}(g_{t})\over \lambda_{k}(g_{t})}\big|_{t=0^-}
=\left[\lambda_{k}(g) {d \over dt}\lambda_{k+1}(g_t)\big|_{t=0^-}- \lambda_{k+1}(g) {d \over dt}\lambda_{k}(g_t)\big|_{t=0^-}\right]$$ 
is the greatest eigenvalue of $\tilde{P}_{k,h}$ on $E_{k}(g) \otimes E_{k+1}(g)$, and  $$\lambda_{k}(g)^{2} {d \over dt}{\lambda_{k+1}(g_{t})\over \lambda_{k}(g_{t})}\big|_{t=0^+}
=\left[\lambda_{k}(g) {d \over dt}\lambda_{k+1}(g_t)\big|_{t=0^+}- \lambda_{k+1}(g) {d \over dt}\lambda_{k}(g_t)\big|_{t=0^+}\right]$$ is the least eigenvalue of $\tilde{P}_{k,h}$ on $E_{k}(g) \otimes E_{k+1}(g)$. Hence, the criticality of $g$ for $\lambda_{k+1}/\lambda_{k}$ is equivalent to the fact that $\tilde{P}_{k,h}$ admits eigenvalues of both signs, which is equivalent to the indefiniteness of $\tilde{Q}_{k,h}.$
\end{proof}
\begin{lemma}\label{lem 51}
Let $g$ be a Riemannian metric on $M$. The two following conditions are equivalent:
\begin{enumerate}
	\item [i)] $\forall h \in S^{2}(M)$, the quadratic form $\tilde{Q}_{k,h}$ is indefinite on $E_{k}(g) \otimes E_{k+1}(g)$.
	\item[ii)] There exist two finite families $\left\{u_{1}, \cdots, u_{p}\right\}\subset E_{k}(g)$ and $\left\{v_{1}, \cdots, v_{q}\right\}\subset E_{k+1}(g)$ of eigenfunctions associated with $\lambda_{k}(g)$ and $\lambda_{k+1}(g)$ respectively, such that
	$$ \sum _{i \leq p} (du_{i} \otimes du_{i} + {1 \over 4} \Delta_{g}u_{i}^{2}\,g)=\sum_ {j \leq q} (dv_{j} \otimes dv_{j} + {1 \over 4} \Delta_{g}v_{j}^{2}\,g).$$
\end{enumerate}
\end{lemma}
\begin{proof}
``(i) implies (ii)'':  Let us introduce  the two following convex cones
$$K_{1}=\left\{\sum_{i \in I}(du_{i} \otimes du_{i}+ {1 \over 4}\Delta_{g}u_{i}^{2}\, g); \, u_{i}\in E_{k}(g), \, I\subset \N, \, I \, \hbox {finite}\right\}\subset S^{2}(M)$$
and
$$K_{2}=\left\{\sum_{j \in J}(dv_{j} \otimes dv_{j}+ {1 \over 4}\Delta_{g}v_{j}^{2}\, g); \, v_{j}\in E_{k+1}(g), \, J\subset \N, \, J \, \hbox {finite}\right\}\subset S^{2}(M).$$
It suffices to prove that $K_{1}$ and $ K_{2}$ have a nontrivial intersection. Indeed, otherwise, applying classical separation theorems, we show the existence of a 2-tensor $h \in S^{2}(M)$ such that, $\forall \ T_{1}\in K_{1},\, T_{1}\ne 0$,
$$ \int_{M}(T_{1},h) v_{g} >0$$
and, $\forall \ T_{2}\in K_{2}$,
$$ \int_{M}(T_{2},h) v_{g}  \leq 0.$$
Therefore, $\forall u \in E_{k}(g)$ and $\forall v \in E_{k+1}(g)$, with $u \ne 0$ and $v \ne 0$, one has $Q_{h}(u)<0$, $Q_{h}(v) \ge 0$ and 
\begin{eqnarray*}
\tilde{Q}_{k,h}(u \otimes v)&=&\lambda_{k+1}(g)\left\|v\right\|^{2}_{L^{2}(g)}Q_{h}(u)-\lambda_{k}(g)\left\|u\right\|^{2}_{L^{2}(g)}Q_{h}(v)\\
&\le& \lambda_{k+1}(g)\left\|v\right\|^{2}_{L^{2}(g)}Q_{h}(u)<0,
\end{eqnarray*}
which implies that $\tilde{Q}_{k,h}$ is negative definite on $E_{k}(g)\otimes E_{k+1}(g)$.\\

\noindent``(ii) implies (i)'':  Let $\left\{u_{i}\right\}_{i \leq p}$ and $ \left\{v_{j}\right\}_{j \leq q}$ be as in (ii). From the identity in (ii), we get, after taking the trace and integrating,
$$ \sum_{i \leq p} \int_{M}|du_{i}|^{2} v_{g}=\sum_{j \leq q} \int_{M}|dv_{j}|^{2} v_{g},$$
which gives,
$$ \lambda_{k}(g) \sum_{i \leq p}\left\|u_{i}\right\|^{2}_{L^{2}(g)}=\lambda_{k+1}(g) \sum_{j \leq q}\left\|v_{j}\right\|^{2}_{L^{2}(g)}.$$
Now,
$$\sum_{i,j}\tilde{Q}_{k,h}(u_{i}\otimes v_{j})=\sum_{i,j}\lambda_{k+1}(g) \left\|v_{j}\right\|^{2}_{L^{2}(g)} Q_{h}(u_{i})- \lambda_{k}(g) \left\|u_{i}\right\|^{2}_{L^{2}(g)} Q_{h}(v_{j}).$$
Assumption (ii) implies that $ \sum_{i \leq p}Q_{h}(u_{i})= \sum_{j \leq q}Q_{h}(v_{j}).$
Therefore,
$$ \sum_{i,j}\tilde{Q}_{k,h}(u_{i}\otimes v_{j})=\left(\sum_{j \leq q}\lambda_{k+1}(g) \left\|v_{j}\right\|^{2}_{L^{2}(g)} - \sum_{i \leq p}\lambda_{k}(g) \left\|u_{i}\right\|^{2}_{L^{2}(g)}\right)\sum_{i \leq p} Q_{h}(u_{i})=0.$$
Hence, $\tilde{Q}_{k,h}$ is indefinite on $E_{k}(g) \otimes E_{k+1}(g)$.
\end{proof}

Consequently, we have proved the following
\begin{theorem}\label{th 51}
A metric $g$ on $M$ is critical for the functional $\displaystyle{{\lambda_{k+1}\over \lambda_{k}}}$ if and only if there exist two families $\left\{u_{1}, \cdots, u_{p}\right\}\subset E_{k}(g)$ and $\left\{v_{1}, \cdots, v_{q}\right\}\subset E_{k+1}(g)$ of eigenfunctions associated with $\lambda_{k}(g)$ and $\lambda_{k+1}(g)$, respectively, such that
\begin{equation}\label{7}
\sum _{i \leq p} du_{i} \otimes du_{i} -\sum_ {j \leq q} dv_{j} \otimes dv_{j} =\alpha g
\end{equation}
for some $\alpha \in \mathcal{C}^{\infty}(M)$, and 
\begin{equation}\label{8}
\lambda_{k}(g) \sum_{i \leq p}u_{i}^{2}-\lambda_{k+1}(g) \sum_{j \leq q}v_{j}^{2}={n-2 \over n} (\sum_{i \leq p}|du_{i}|^{2}-\sum_{j \leq q}|dv_{j}|^{2}).
\end{equation}
\end{theorem}
Indeed, a straightforward calculation shows that the two equations (\ref{7}) and (\ref{8}) are equivalent to the condition (ii) of Lemma \ref{lem 51}.
\begin{corollary}
If $g$ is a critical metric of the functional ${\lambda_{k+1} \over \lambda_{k}}$, with $\lambda_{k+1}(g) \ne \lambda_{k}(g)$, then 
   $$ \inf \left\{\dim E_{k}(g),\dim E_{k+1}(g)\right\} \geq 2. $$
\end{corollary}
\begin{proof}
Let $\left\{u_{i}\right\}_{i \leq p}\subset E_{k}(g)$ and $\left\{v_{j}\right\}_{j \leq q}\subset E_{k+1}(g)$
be two families of eigenfunctions satisfying (\ref{7}) and (\ref{8}) above. Taking the trace in (\ref{7}) and using (\ref{8}) we get 
\begin{equation}\label{9}
\alpha= {1 \over n}(\sum_{i \leq p}|du_{i}|^{2}-\sum_{j \leq q}|dv_{j}|^{2})={1 \over 4}\Delta_{g}(\sum_{j \leq q}v_{j}^{2}- \sum_{i \leq p}u_{i}^{2}).
\end{equation}

Assume that $\alpha=0$. Using (\ref{8}) and (\ref{9}), we deduce that both $\sum_{i \leq p}u_{i}^{2}$ and $\sum_{j \leq q}v_{j}^{2}$ are constant on $M$. Since the $u_{i}$'s and the $v_{j}$'s are not constant, we get the result.
Assume now $\alpha \ne 0$. Since $\int_{M} \alpha\ v_{g}=0$ (see (\ref {9})), the function $\alpha$ takes both positive and negative values. Let $x \in M$ such that $\alpha(x)>0$. From (\ref{7}), the quadratic form $\sum _{i \leq p} du_{i} \otimes du_{i} $ is clearly positive definite on $T_{x}M$. Hence, the family 
$\left\{du_{i}\right\}$ has maximal rank at $x$. This shows that $\dim E_{k}(g) \geq n$. At a point $x \in M$ where $\alpha(x)<0$, the quadratic form  
$\sum_ {j \leq q} dv_{j} \otimes dv_{j}$ is positive definite on $T_{x}M$ and, then, $\dim E_{k+1}(g) \geq n$.
\end{proof}

When we deal with critical metrics of the functional ${\lambda_{k+1} \over \lambda_{k}}$ restricted to $C(g)$, only tensors of the form $h= \varphi\, g$, with $\varphi \in \mathcal{C}^{\infty}(M)$, are involved. The corresponding quadratic forms on $E_{k}(g) \otimes E_{k+1}(g)$ are  given by
$$ \tilde{q}_{k,\varphi}(u \otimes v)= \lambda_{k+1}(g)\left\|v\right\|^{2}_{L^{2}(g)} q_{\varphi}(u)-\lambda_{k}(g)\left\|u\right\|^{2}_{L^{2}(g)} q_{\varphi}(v).$$
Following the steps of the proof of Proposition \ref{prop 51}, we can show that:
\begin{proposition}
A Riemannian metric $g$ on $M$ is critical for the functional ${\lambda_{k+1} \over \lambda_{k}}$ restricted to $C(g)$ if and only if, $\forall \varphi \in \mathcal{C}^{\infty}(M)$, the quadratic form $\tilde{q}_{k,\varphi}$ is indefinite on $E_{k}(g) \otimes E_{k+1}(g)$.
\end{proposition}

Replacing the convex cones $K_{1}$ and $K_{2}$ in the proof of Lemma \ref{lem 51} by
$$ H_{1}= \left\{\sum_{i \in I}(\lambda_{k}(g) u_{i}^{2}- {n-2 \over n} |du_{i}|^{2}); \, u_{i}\in E_{k}(g),\,I \subset \N,\,  I \hbox{ finite}\right\}\subset L^{2}(M,g)$$
and
$$ H_{2}= \left\{\sum_{j \in J}(\lambda_{k+1}(g) v_{j}^{2}- {n-2 \over n} |dv_{j}|^{2}); \, v_{j}\in E_{k+1}(g),\,J \subset \N,\,  J \hbox{ finite}\right\}\subset L^{2}(M,g),$$
we can show, by the same arguments, that the indefiniteness of $\tilde{q}_{k,\varphi}$ for all $\varphi \in \mathcal{C}^{\infty}(M)$, is equivalent to the fact that $H_{1}$ and $H_{2}$ have a non-trivial intersection. Therefore, one has:

\begin{theorem}\label{th 52}
A Riemannian metric $g$ on $M$ is critical for the functional ${\lambda_{k+1} \over \lambda_{k}}$ restricted to $C(g)$ if and only if, there exist two families 
$\left\{u_{1}, \cdots, u_{p}\right\}\subset E_{k}(g)$ and $\left\{v_{1}, \cdots, v_{q}\right\}\subset E_{k+1}(g)$ of eigenfunctions associated with $\lambda_{k}(g)$ and $\lambda_{k+1}(g)$, respectively, such that
$$\lambda_{k}(g) \sum_{i \leq p}u_{i}^{2}-\lambda_{k+1}(g) \sum_{j \leq q}v_{j}^{2}={n-2 \over n} (\sum_{i \leq p}|du_{i}|^{2}-\sum_{j \leq q}|dv_{j}|^{2}).$$
\end{theorem}

\begin{remark}
In dimension $2$, the condition of Theorem \ref{th 52} amounts to 
$$ \sum_{j \leq q}v_{j}^{2}= {\lambda_{k}(g) \over \lambda_{k+1}(g)}\sum_{i \leq p} u_{i}^{2}.$$It is clear that in this case, if $\lambda_{k+1}(g) \ne \lambda_{k}(g)$, then at least one of the eigenvalues $\lambda_{k}(g)$ and $\lambda_{k+1}(g)$ is degenerate.
\end{remark}

\end{document}